\documentclass[conference]{IEEEtran}
\IEEEoverridecommandlockouts

\usepackage{cite}
\usepackage{amsmath,amssymb,amsfonts}
\usepackage{algorithmic}
\usepackage{graphicx}
\usepackage{textcomp}
\usepackage[dvipsnames]{xcolor}
\def\BibTeX{{\rm B\kern-.05em{\sc i\kern-.025em b}\kern-.08em
    T\kern-.1667em\lower.7ex\hbox{E}\kern-.125emX}}
\begin{document}

\newtheorem{theorem}{Theorem}
\newtheorem{corollary}{Corollary}

\newenvironment{reptheorem}[1]
  {\renewcommand\therthm{\ref*{#1}}\rthm}
  {\endrthm}

%
%

\newcommand{\veca}{\boldsymbol{a}}
\newcommand{\vecb}{\boldsymbol{b}}
\newcommand{\vecc}{\boldsymbol{c}}
\newcommand{\vecd}{\boldsymbol{d}}
\newcommand{\vece}{\boldsymbol{e}}
\newcommand{\vecf}{\boldsymbol{f}}
\newcommand{\vecg}{\boldsymbol{g}}
\newcommand{\vech}{\boldsymbol{h}}
\newcommand{\veci}{\boldsymbol{i}}
\newcommand{\vecj}{\boldsymbol{j}}
\newcommand{\veck}{\boldsymbol{k}}
\newcommand{\vecl}{\boldsymbol{l}}
\newcommand{\vecm}{\boldsymbol{m}}
\newcommand{\vecn}{\boldsymbol{n}}
\newcommand{\veco}{\boldsymbol{o}}
\newcommand{\vecp}{\boldsymbol{p}}
\newcommand{\vecq}{\boldsymbol{q}}
\newcommand{\vecr}{\boldsymbol{r}}
\newcommand{\vecs}{\boldsymbol{s}}
\newcommand{\vect}{\boldsymbol{t}}
\newcommand{\vecu}{\boldsymbol{u}}
\newcommand{\vecv}{\boldsymbol{v}}
\newcommand{\vecw}{\boldsymbol{w}}
\newcommand{\vecx}{\boldsymbol{x}}
\newcommand{\vecy}{\boldsymbol{y}}
\newcommand{\vecz}{\boldsymbol{z}}
\newcommand{\vecpi}{\boldsymbol{\pi}}
\newcommand{\vecxi}{\boldsymbol{\xi}}
\newcommand{\veczero}{\boldsymbol{0}}
\newcommand{\vecone}{\boldsymbol{1}}
\newcommand{\vecbeta}{\boldsymbol{\beta}}
\newcommand{\vectheta}{\boldsymbol{\theta}}


\newcommand{\matrixW}{\mathcal{W}}
\newcommand{\matrixT}{\mathcal{T}}
\newcommand{\matrixM}{\mathcal{M}}
\newcommand{\matrixD}{\mathcal{D}}

\newcommand{\barvecx}{\Bar{\boldsymbol{x}}}
\newcommand{\bartheta}{\Bar{\theta}}
\newcommand{\barvecy}{\Bar{\boldsymbol{y}}}

\newcommand{\setA}{\textit{A}}
\newcommand{\setN}{\textit{N}}
\newcommand{\setS}{\textit{S}}
\newcommand{\setV}{\textit{V}}
\newcommand{\setX}{\textit{X}}
\newcommand{\setY}{\textit{Y}}
\newcommand{\setGamma}{\varGamma}
\newcommand{\arcSet}{\textit{A}}
\newcommand{\nodeSet}{\textit{N}}
\newcommand{\setPi}{\varPi}
\newcommand{\setvpi}{\textit{V}(\setPi)}
\newcommand{\hatsetV}{\hat{\setV}}
\newcommand{\setOmega}{\varOmega}

\newcommand{\sublevelset}{\leftindex_{s}{\mbox{S}}}
\newcommand{\levelset}{\mbox{S}}
\newcommand{\levelEV}{\levelset_{EV}}
\newcommand{\levelBM}{\levelset_{BM(\hatV)}}
\newcommand{\levelPV}{\levelset_{PV}}
\newcommand{\levelEF}{\levelset_{EF}}

\newcommand{\AOSKernel}{\mbox{AOSKernel}}
\newcommand{\epigraph}{\mbox{epi}}
\newcommand{\conv}{\mbox{conv}}
\newcommand{\proj}{\mbox{proj}}
\newcommand{\projX}{\proj_{\mathR^{n_1}}}
\newcommand{\funcFS}{\mbox{FS}}
\newcommand{\funcRS}{\mbox{RS}}

\newcommand{\funca}{\mbox{a}}
\newcommand{\funcA}{\mbox{A}}
\newcommand{\funcb}{\mbox{b}}
\newcommand{\funcB}{\mbox{B}}
\newcommand{\funcc}{\mbox{c}}
\newcommand{\funcC}{\mbox{C}}
\newcommand{\funcd}{\mbox{d}}
\newcommand{\funcD}{\mbox{D}}
\newcommand{\funce}{\mbox{e}}
\newcommand{\funcE}{\mbox{E}}
\newcommand{\funcf}{\mbox{f}}
\newcommand{\funcF}{\mbox{F}}
\newcommand{\funcg}{\mbox{g}}
\newcommand{\funcG}{\mbox{G}}
\newcommand{\funch}{\mbox{h}}
\newcommand{\funcH}{\mbox{H}}
\newcommand{\funci}{\mbox{i}}
\newcommand{\funcI}{\mbox{I}}
\newcommand{\funcj}{\mbox{j}}
\newcommand{\funcJ}{\mbox{J}}
\newcommand{\funck}{\mbox{k}}
\newcommand{\funcK}{\mbox{K}}
\newcommand{\funcl}{\mbox{l}}
\newcommand{\funcL}{\mbox{L}}
\newcommand{\funcm}{\mbox{m}}
\newcommand{\funcM}{\mbox{M}}
\newcommand{\funcn}{\mbox{n}}
\newcommand{\funcN}{\mbox{N}}
\newcommand{\funco}{\mbox{o}}
\newcommand{\funcO}{\mbox{O}}
\newcommand{\funcp}{\mbox{p}}
\newcommand{\funcP}{\mbox{P}}
\newcommand{\funcq}{\mbox{q}}
\newcommand{\funcQ}{\mbox{Q}}
\newcommand{\funcr}{\mbox{r}}
\newcommand{\funcR}{\mbox{R}}
\newcommand{\funcs}{\mbox{s}}
\newcommand{\funcS}{\mbox{S}}
\newcommand{\funct}{\mbox{t}}
\newcommand{\funcT}{\mbox{T}}
\newcommand{\funcu}{\mbox{u}}
\newcommand{\funcU}{\mbox{U}}
\newcommand{\funcv}{\mbox{v}}
\newcommand{\funcV}{\mbox{V}}
\newcommand{\funcw}{\mbox{w}}
\newcommand{\funcW}{\mbox{W}}
\newcommand{\funcx}{\mbox{x}}
\newcommand{\funcX}{\mbox{X}}
\newcommand{\funcy}{\mbox{y}}
\newcommand{\funcY}{\mbox{Y}}
\newcommand{\funcz}{\mbox{z}}
\newcommand{\funcZ}{\mbox{Z}}

\newcommand{\funchatQV}{\hatQVVariant{\hatV}}

\newcommand{\hatV}{\hat{\setV}}
\newcommand{\hatVT}{\hatV_t}
\newcommand{\vpi}{\funcV(\setPi)}
\newcommand{\mathR}{\mathbb{R}}
\newcommand{\hatQVVariant}[1]{\funcQ_{#1}}
\newcommand{\hatQV}{\hatQVVariant{\hatV}}
\newcommand{\hatQVT}{\hatQVVariant{\hatV_t}}
\newcommand{\hattheta}{\hat{\theta}}
\newcommand{\hatvecy}{\hat{\vecy}}
\newcommand{\hatsetY}{\hat{\setY}}

\newcommand{\defeq}{:=}
\newcommand{\BMVariant}[1]{BM(#1)}
\newcommand{\BMhatV}{\BMVariant{\hatV}}
\newcommand{\BMVT}{\BMVariant{\hatV_t}}

\newcommand{\labelspace}{\mbox{ }\quad\mbox{ }\quad\mbox{ }}

\newcommand{\aosExactSetName}{EX-ALT}
\newcommand{\aosApproxSetName}{A-ALT}
\newcommand{\nameEF}{extensive-form}
\newcommand{\namePV}{projected\ variable}
\newcommand{\nameEV}{epigraphical\ variant}
\newcommand{\nameBM}{Benders\ master}
\newcommand{\nameBenders}{Benders\ Decomposition}
\newcommand{\namePrimalValueFunction}{primal\ value\ function}
\newcommand{\nameDualValueFunction}{dual\ value\ function}
\newcommand{\DCI}{Brown et al.}
\newcommand{\IW}{Israeli-Wood}
\newcommand{\nameMaster}{master}
\newcommand{\nameEnumerateLinear}{Pyomo-AOS-Linear}
\newcommand{\nameEnumerateBinary}{Pyomo-AOS-Binary}

\newcommand{\fillInCitation}{[FIC]}
\newcommand{\fillInDataPoint}{XYZ}

\newcommand{\ppp}{policy planning problem}
\newcommand{\ppps}{policy planning problems}
\newcommand{\Ppp}{Policy planning problem}
\newcommand{\Ppps}{Policy planning problems}
\newcommand{\PPP}{Policy Planning Problem}
\newcommand{\PPPs}{Policy Planning Problems}

\newcommand{\altsol}{alternative solution}
\newcommand{\Altsol}{Alternative solution}
\newcommand{\altsols}{alternative solutions}
\newcommand{\Altsols}{Alternative solutions}
\newcommand{\equationRef}[1]{(\ref{#1})}

\newcommand{\JKS}[1]{\textcolor{MidnightBlue}{\textbf{[JKS: #1]}}}
\newcommand{\MPV}[1]{\textcolor{YellowOrange}{\textbf{[MPV: #1]}}}

\newcommand{\TODO}[1]{\textbf{TODO: #1} \\}
%
%

\title{An Optimal Solution is Not Enough: Alternative Solutions and Optimal Power Systems
\thanks{This work was supported in part by the Laboratory Directed Research and Development program at Sandia
National Laboratories, a multimission laboratory managed and operated by National Technology and
Engineering Solutions of Sandia LLC, a wholly owned subsidiary of Honeywell International Inc. for the U.S.
Department of Energy’s National Nuclear Security Administration under contract DE-NA0003525.	Additional funding came from the Morgridge Chair of Computer Sciences at University of Wisconsin-Madison. 

979-8-3315-4112-5/25/\$31.00 ©2026 IEEE}
}

\author{\IEEEauthorblockN{Matthew P. Viens}
\IEEEauthorblockA{\textit{Center for Computing Research} \\
\textit{Sandia National Laboratories}\\
Albuquerque, NM, USA \\
mpviens@sandia.gov}
\and
\IEEEauthorblockN{J. Kyle Skolfield and William E. Hart}
\IEEEauthorblockA{\textit{Center for Computing Research} \\
\textit{Sandia National Laboratories}\\
Albuquerque, NM, USA \\
\{jkskolf,wehart\}@sandia.gov}
\and
\IEEEauthorblockN{Michael C. Ferris}
\IEEEauthorblockA{\textit{Department of Computer Sciences} \\
\textit{University of Wisconsin-Madison}\\
Madison, WI, USA \\
ferris@cs.wisc.edu}
}
\maketitle

\begin{abstract}
Power systems modeling and planning has long leveraged mathematical programming for its ability to provide optimality and feasibility guarantees. One feature that has been recognized in the optimization literature since the 1970s is the existence and meaning of multiple exact optimal and near-optimal solutions, which we call alternative solutions. In power systems modeling, the use of alternative solutions has been limited to energy system optimization modeling (ESOM) applications and modeling to generate alternative (MGA) techniques. We present three key results about alternative solutions for power systems modeling. First, we give a perspective, based on sublevel sets and projection, for characterizing alternative solutions as a facet of general optimization theory. Second, we include pointers to alternative solution generation methods and tools beyond MGA-style techniques. Third, we demonstrate the use cases for alternative solutions in power system modeling on the fundamental optimal power flow problem.

\end{abstract}


\section{Introduction}
Power systems has long been a core application area for computational modeling and computational decision-making. Mathematical optimization has been essential to discovering minimal-cost actions and other best-objective style decisions. One recurring task in the intersection between mathematical optimization and power systems is the need to continually bring in cutting-edge advancements in computational solvers, theoretical results, and modeling tools to address the ever-expanding demands on our power-system resources in the electricity grid. Our work extends this tradition by presenting alternative optimal solutions as an essential part of optimal power systems modeling.

Alternative (optimal) solutions (AOS) is our term for the idea that both decision-making problems and mathematical optimization models can have multiple solutions that achieve an optimal or near-optimal objective result. The mathematical optimization literature recognized the value of generating such alternative solutions in the 1970s via separate lines of research in infrastructure planning and agricultural economics, from Brill in \cite{brill_1979_public_planning_alternatives} and Paris in \cite{paris_mult_opt_lp}, respectively. From the advent of better optimization solvers and cheaper computational resources, there has been a resurgence of interest in generating alternative solutions for the last 15 years, notably in the planning of future electrical grids in \cite{VOLL_optimum_is_not_enough} and in \cite{decarolis_alternatives_mga_energy_futures}. AOS methods for design of future electrical grids have largely been focused on capacity expansion (CapEx) and energy-system modeling (ESOM) especially in the context of integrating renewables-based generation. The existing ESOM communities that leverage AOS generation methods derive from Brill's line of work and inherit his terminology of modeling to generate alternatives (MGA), see reviews in \cite{Lau_2024_mga_alternatives_energy_applications_review} and \cite{turan2025oracle_preprint}. We prefer the AOS terminology over MGA for several reasons, including: there are more uses for alternative solutions than solely generation (e.g. analysis of solution sets), `alternative solutions' is clearer than just `alternatives', and AOS is consistent with phrasing in the textbooks that treat these concepts in \cite{paris_economic_foundations_of_symmetric_programming} and \cite{hpwilliams_model_building}.

In this work, we present three key concepts regarding alternative solutions. In Section \ref{paradigm}, we present a framework for understanding alternative solutions as a component of the larger theory of mathematical optimization. In Section \ref{generation}, we present methodological and software resources for generation and use of alternative solutions to match a range of modeling structures and model use cases. In Section \ref{application}, we demonstrate how alternative solutions can impact a range of optimal power systems problems by focusing on a model core to many applications: optimal power flow. We conclude with a summary of results and discussion of further applications in Section \ref{conclusion}.
\section{Alternative Solutions as a Paradigm}
\label{paradigm}
\subsection{Problem Statement}
\label{paradigm:problem_statement}
When we talk about alternative (optimal) solutions, we first need to recognize that solutions to decision-making problems are not automatically unique. There may be two ways to traverse a network that take the same amount of time or multiple generation expansion plans that meet future electric demand at equal cost. There also may be additional ways to traverse that network that takes slightly more time but use a different traversal path.  Likewise, it is easy to envision a plethora of expansion plans that maintain network reliability at marginally increased cost. When we move to using mathematical models of decision-making, we still have this possibility for multiple solutions, and our decision-making tools should reflect this reality. Alternative optimal solutions is the term we use to to capture this concept in mathematical optimization.

There are many reasons to care about the existence of alternative solutions. One reason is to understand how a particular solution is picked when alternative solutions may exist. Another reason is to address secondary concerns or metrics not included in a model -- for example, capturing the impacts of public opinion -- by evaluating alternative solutions using those metrics.  An additional reason is to use alternative solutions to discover, rank, or relate secondary objectives and tradeoffs in a decision making model. All of these reasons can occur in optimal power systems.
\subsection{What \& So What}
\label{paradigm:what_and_so_what}
For a concrete example of what alternative optimal solutions can look like mathematically, here is a simple optimization problem:
\begin{subequations}
\begin{align}
    \max_{x \in \mathR^2} x_1 \\
    x_2 \leq -x_1 + 101 \\
    \label{eq:vert_line}
    x_1 \leq 100 \\
    x_2 \leq 100
\end{align}
\end{subequations}
The feasible region $X$ of this problem can be seen in Figure 1.

\begin{figure}
    \centering
    \includegraphics[width=0.4\textwidth]{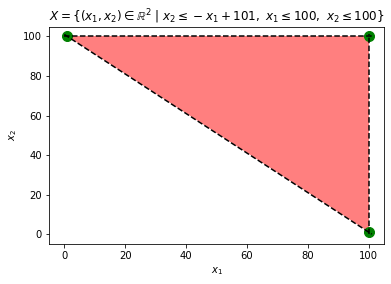}
    \caption{Right-angle triangle used for feasibility region.}
    \label{fig:feas_region_1}
\end{figure}

There is exactly one optimal objective value: 100. The set of points that achieve this optimal objective value is 
$S_1 = \{x | x_1 = 100, x_2 \in [1,100]\}$; this set clearly contains an infinite number of points. 

This simple example points to two key ideas: infinite alternative solutions and modeling limits. For the case where infinite alternative solutions exist, there are techniques to create finite representations, which is addressed more in Section \ref{paradigm:analysis:representations}. For modeling limits, there are often clear bounds on what can be included in optimization models both in terms of additional constraints and secondary concerns. As an example, consider the secondary concern of maximizing $x_2$. The best possible result for both of these concerns is then $x_1 = 100, x_2 = 100$, achieving the best score of 100 on both the first and second objective. The worst possible objective for this secondary concern is at $x_1 = 100, x_2 = 1$. This means that both the best and worst answers for secondary concerns may be in the set of points that achieve the optimal first objective value with a solver arbitrarily picking from $S_1$. 

There are a few ways of addressing this secondary objective concern. First, if the secondary condition is known when solving the original problem and can be treated by optimization tools, there are techniques like multi-objective (or Pareto) optimization. Second, alternative solutions from the $S_1$ set can be identified using the techniques described in Section \ref{generation}. The second case covers circumstances in which the secondary concerns are not known a priori or cannot be included in an optimization model (e.g. ranking by public utility commission).

In the previous example, we were dealing with exact alternative solutions. There are also use cases for which only considering alternative solutions that exactly achieve the optimal objective value is quite restrictive. This can be seen in the following problem:
\begin{subequations}
\begin{align}
    \max_{x \in \mathR^2} x_1 \\
    x_2 \leq -x_1 + 101 \\
    \label{eq:near_vert_line}
    x_2 \leq -99x_1 + 9901 \\
    x_2 \leq 100
\end{align}
\end{subequations}
This results in a slightly different feasible region as see in Figure 2.

\begin{figure}
    \centering
    \includegraphics[width=0.4\textwidth]{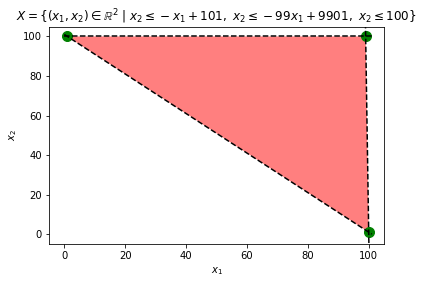}
    \caption{Slightly obtuse triangle used for feasibility region.}
    \label{fig:feas_region_2}
\end{figure}

Figure 2 represents an extremely slight shift from Figure 1: there is a slight perturbation in the location of the upper right vertex from $(100,100)$ to $(99,100)$ due to the change from Equation \ref{eq:vert_line} to Equation \ref{eq:near_vert_line}. Small differences in constraints can arise from different data-sources used, different estimates of costs, or many other factors. Even this slight change in data markedly changes the optimal outcomes of our problem. The set of points that achieves the best objective result of 100 on the first objective is now a singleton set $S_2 = {(100,1)}$; the impact of perturbations on alternative solutions is one element of stability theory \cite{robinson_stability_system_of_inequalities} \cite{robinson_stability_linear_programs}.

When the same secondary concern, maximizing $x_2$, is applied to this perturbed problem, the changes are even more pronounced.  If only solutions in $S_2$ are used -- i.e., those that exactly achieve the best answer on the primary objective -- there is only one option, which achieves a secondary concern value of $1$. When near-optimal solutions, even just $1\%$ from optimal, are allowed, then a 100x improvement on the secondary objective can be achieved with $x_1 = 99, x = 100$. This demonstrates that a general definition of alternative solutions should be capable of including near-optimal solutions.
\subsection{Analysis of Alternative Solutions}
\label{paradigm:analysis}
Given a general optimization problem
\begin{align*}
    z^* = \min_{x \in X} f(x),
\end{align*}
for which there exists an optimal solution, there are three key tools to analyze and understand alternative solutions: sublevel sets, minimal representations, and projections.
\subsubsection{Sublevel Sets}
\label{paradigm:analysis:sublevel_sets}
There are several ways to describe sets of feasible points and feasible points that achieve a specific objective value. We standardize on the following sublevel set framework:
\begin{align*}
    \levelset (f,X,\tau) &= \{x \in X\ | \ f(x) \leq \tau\}.
\end{align*}
This defines the set of all points in the feasible space that achieve a function value less than $\tau$, which we call the level value. When $\tau < z^*$, the result is the empty set. When $\tau = z^*$, the resulting set is non-empty and takes the form:
\begin{align*}
    \levelset (f,X,z^*) &= \{x \in X\ | \ f(x) \leq z^*\} \\
    &= \{x \in X\ | \ f(x) = z^*\} \cup \{x \in X\ | \ f(x) < z^*\} \\
    &= \{x \in X\ | \ f(x) = z^*\}.
\end{align*}
As a result, this general structure allows description of both exact and near-optimal alternative solutions depending on the value of $\tau$. This sublevel set definition is implicitly used but not named in many papers for AOS generation methods with description of points as $f(x) \leq z^* + \epsilon$ \cite{brill_1982_alternatives_HSJ_hop_skip_jump}, \cite{paris_mult_opt_lp}, \cite{decarolis_alternatives_mga_energy_futures}, \cite{Lau_2024_mga_alternatives_energy_applications_review}. Recognizing and naming this critical object allows usage of the general theory that the mathematical optimization literature has developed around sublevel sets, including Rockafellar's advanced treatment \cite{Convex_Analysis_Rockafellar} and Robinson's stability treatment \cite{robinson1984_stability_sublevel_set}. One important case is when sublevel sets are convex, a property called quasiconvexity. There are large sets of functions known to be quasiconvex including trivial examples like linear and convex functions and surprising examples like the square-root and floor functions \cite{royset2021optimization_primer}.

\subsubsection{Minimal Representations}
\label{paradigm:analysis:representations}
The representation of alternative solutions gives a general way of analyzing such points, but does not necessarily give an easily constructible or computationally amenable structure. There are two important cases of sublevel sets to address: singleton sets, and efficient representations. 

For the case of singleton sets, this corresponds to two cases where there will be a unique minimizer. When $\tau = z^*$ and the sublevel set is a singleton, the optimization model has a unique solution. When $\tau > z^*$ and $\levelset(f,X,\tau)$ is a singleton, not only is there a unique solution to the optimization model but there are no other feasible points that achieve an objective value between $z^*$ and $\tau$. This case is much more common in optimization problems involving discrete variables. A natural question would be: what do these singleton sets have to do with alternative solutions? An AOS generation method that is exhaustive and run with level $\tau$ returning only a single point gives a certificate of the sublevel set being a singleton and the uniqueness of the minimizer. This is a broadly applicable test for uniqueness for problem classes that do not have a uniqueness guarantee (e.g., integer programs) compared to those that do (e.g., strongly convex programs).

For efficient representations, there are classes of optimization problems that have finite representations of the feasible region and therefore the sublevel sets. For linear programs (LPs), the feasible region, and sublevel sets, are polyhedral. This means there is a finite representation in terms of the extreme points and extreme rays (as a result of the Minkowski-Weyl theorem \cite{convex_opt_theory}). For mixed-integer linear programs (MILPs), there is an equivalent representation with a specific relaxation called a perfect formulation \cite{conforti2014integer}. In both cases, there is a compact way to represent the sublevel set for a given $\tau$ value and therefore the associated alternative solutions. For other classes the theory is sparser.  For example, quadratic programming (QP) has some results specifically for $\tau = z^*$ \cite{paris_economic_foundations_of_symmetric_programming}. There are AOS generation methods that take advantage of this finite representation \cite{lee2000recursive} and those that do not \cite{brill_1982_alternatives_HSJ_hop_skip_jump}. For the methods that leverage the finite representation, AOS generation can be reduced to either enumerating this representation, if the number of points is small enough, or sampling this representation, in larger cases. The enumeration case is demonstrated in Section \ref{application:three_bus}.

\subsubsection{Projection}
\label{paradigm:analysis:projection}
When there is an interesting subset of variables, then a natural question is: what happens when projecting onto this subset? We define a projection operator onto the subspace of the first $k$ variables:
\begin{equation}
    \label{eqn:proj_x_def}
    \proj_k(\vecp) \defeq \matrixM \vecp, \matrixM \in \mathR^{k \times \dim(\vecp)},\matrixM_{ij} = \left\{\begin{matrix}
                      1 \quad i=j \\
                      0 \quad i\neq j
                    \end{matrix} \right .
\end{equation}
This projection operation is defined for points and can be applied point-wise to the sublevel set representation of alternative solutions.

In the case where the majority of the decision-making meaning of the solution is contained in a subset of the variables, projecting alternative solutions onto that subset can be useful both for theoretical analysis and for actual generation of alternative solutions. Both use cases are demonstrated in the example in Section \ref{application}, focusing on generation variables in a DC optimal power flow problem (DC-OPF). Generating alternative solutions in the non-projected space and then projecting it leads to challenges in computation and potential redundancies. However, there are methods for specific problem structures (e.g. stochastic programming) that can perform AOS generation directly in the projected space \cite{aos_benders_preprint_viens_hart_ferris_2025}. 
\section{Generation of Alternative Solutions}
\label{generation}
We split our discussion of generation of alternative solutions into theoretical methods and software approaches.
\subsection{Theoretical Methods}
\label{generation:methods}
Although the details of the many AOS generation theoretical methods is beyond the scope of this paper, it is important to present a brief overview to identify common trends. The literature on AOS generation methods is split amongst several research communities that independently developed techniques. To the best of the authors' knowledge, there is not a comprehensive review of methods either generally or by problem class. We provide pointers to comprehensive reviews by domain, highlight common techniques, and point to techniques not otherwise emphasized in the literature. 

There are review of techniques from the ESOM literature using what they call MGA techniques. In the context of reviewing performance on ESOM at scale, there is a review from 2024 in \cite{Lau_2024_mga_alternatives_energy_applications_review}. In the context of proposing a new method, there is a review from 2025 in \cite{turan2025oracle_preprint}. The ESOM use case generally focuses on LPs, but the random vector and variable max techniques can be applied generally. Additionally, the term MGA is derived from early work by Brill for infrastructure problems in with a method called Hop-Skip-Jump \cite{brill_1982_alternatives_HSJ_hop_skip_jump}. Another source of early work was Paris for agricultural economics applications, his techniques for LPs, QPs, and complementarity problems are given a combined treatment as a textbook chapter in \cite{paris_economic_foundations_of_symmetric_programming}. 

For MILPs and combinatorial problems, there is a review of generation methods in the context of `solution engineering' in \cite{Petit2019_alternative_solutions_as_solution_engineering_review_ijoc}. Diversification of solutions in the MILP context is treated in \cite{Ahanor_Trapp_DiversiTree_IJOC} and in \cite{danna2007_multiple_solutions_for_mixed_integer_programs_mips_ILOG_CPLEX} with the latter providing the classic `one-tree' algorithm.  

For further methods, there are techniques for LPs that generate the extreme point representation in \cite{lee2000recursive}. There are methods for binary and integer programs that rely on the general concept of no-good cuts in \cite{balas_no_good_cuts} and in \cite{denegre2011interdiction_integer_no_good_cuts_thesis}. There are methods that generate under projections \cite{aos_benders_preprint_viens_hart_ferris_2025}, which can be useful for stochastic programming and bilevel programming. There are also solver-specific techniques that are not publicly known, which are described more in the software discussion below.
\subsection{Software Approaches}
\label{generation:software}
While there are many theoretical approaches to generation of alternative solutions, there are concrete software generation tools. We focus here on supported software methods that are part of either active open-source projects or commercial tools. There are unsupported software projects that may be published as part of the journal publications for generation methods described in the previous section.

For supported software, there are several classes. For application-specific software, the PyPSA family of ESOM models have several tools either in PyPSA core or in related models like PyPSA-eu \cite{PyPSA}. For optimization solvers, both Gurobi \cite{gurobi} and CPLEX \cite{cplex} have the ability to generate alternative solutions under construct called `solution pool'. Since both are commercial solvers, their methods are functionally black-box. As both are primarily MILP solvers, there is some nuance and ambiguity in how differences between alternative solutions over the continuous variables are handled. For optimization modeling languages, both GAMS \cite{gams} and Pyomo \cite{bynum2021pyomo_book} have the ability to connect to solution pool capabilities in a specific solver. Pyomo also has implementations of both LP extreme point discovery and binary program no-good cuts methods that are agnostic to the underlying solver \cite{cimor_aos_report}.
\section{Application: Optimal Power Flow}
\label{application}
We consider the following DC Optimal Power Flow (DC-OPF) problem \cite{frank2012optimal_power_flow_standard_reference_opf} on Graph $G=(N, E)$:
\begin{subequations}
\label{eq:dc_opf:full}
\begin{align}
\label{eq:dc_opf:objective}
z^* = &\min_{P,f,\theta} \sum_{i \in N} c_i P_i \\
\label{eq:dc_opf:flow_balance}
P_i + \sum_{(j,i) \in E} f_{j,i} - \sum_{(i,j)} f_{i,j} &= L_i, \forall i \in N \\
\label{eq:dc_opf:angle_constraint}
\frac{1}{x_{i,j}}(\theta_i - \theta_j) &= f_{i,j}, \forall (i,j) \in E \\
\label{eq:dc_opf:gen_domain}
P_i &\in [0, P^{max}_i], \forall i \in N \\
\label{eq:dc_opf:flow_domain}
f_{e} &\in [-f^{max}_e, f^{max}_e], \forall e \in E \\
\label{eq:dc_opf:angles_domain}
\theta_i &\in \mathR, \forall i \in N
\end{align}
\end{subequations}
We write the constraints for DC-OPF compactly as $$X_{DC} = \{(P,f,\theta)\ |\ (\ref{eq:dc_opf:flow_balance}),\ldots,(\ref{eq:dc_opf:angles_domain})\}$$ and the objective as $g(P) = c^T P$. $P$ denotes generation at each node, $f$ the flow between nodes, and$\theta$ the bus voltage angles. The data is $L$ for load/demand. We specify the order of the variables as $(P,f,\theta)$ for use in later projection operators.

The Network Flow problem can be made as a relaxation of the DC-OPF problem with constraints:
$$X_{NF} = \{(P,f)\ |\ (\ref{eq:dc_opf:flow_balance}), (\ref{eq:dc_opf:gen_domain}), (\ref{eq:dc_opf:flow_domain})\}.$$
Similarly, the Copper Plate problem can be made as a further relaxation, but it requires a new constraint to replace the earlier flow balance constraint:
\begin{align}
\label{eq:copper_plate:gen_balance}
    \sum_{i \in N} P_i = \sum_{i \in N} L_i,
\end{align}
then the Copper Plate constraints are:
$$X_{CP} = \{P\ |\ (\ref{eq:dc_opf:gen_domain}), (\ref{eq:copper_plate:gen_balance})\}.$$
Note that DC-OPF, Network Flow, and Copper Plate all have the same objective $g(P)$.

It is clear from this presentation that Network Flow is a relaxation of DC-OPF, and Copper Plate is a relaxation of both. Using the sublevel sets and projection operators developed in Section \ref{paradigm:analysis}, this result can be extended to the alternative solutions. For any level value of the objective $g(P)$, the Network Flow sublevel set contains the DC-OPF sublevel set when projected to the relaxed space; the Copper Plate sublevel set contains both the Network Flow and DC-OPF sublevel sets when similarly projected. Since this result holds for any level value for the objective, it also applies specifically to the exact alternative solutions. The full statement and proofs of these results are in Appendix \ref{theorem_appendix}.
\subsection{3-Bus Optimal Power Flow}
\label{application:three_bus}
We consider the 3-Bus system described in Figure 3.

\begin{figure}
    \centering
    \includegraphics[width=0.4\textwidth]{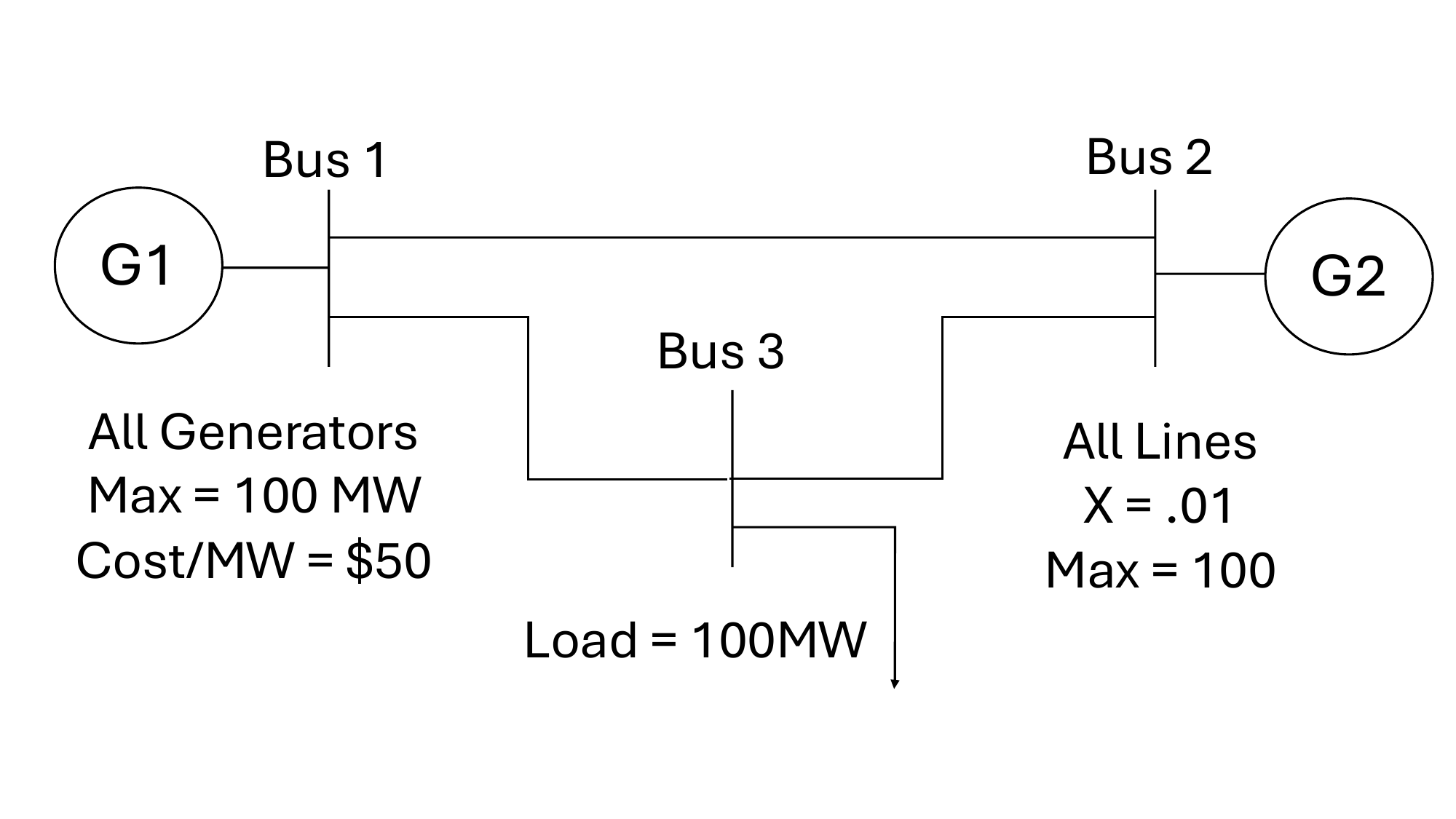}
    \caption{Three-bus symmetric network used for DC-OPF Exemplar.}
    \label{fig:tiny_opf}
    \vspace{-1em}
\end{figure}

Even in this extremely simple optimal power flow problem, alternative solutions will appear. While this is a toy example intended to illustrate the concept, in fact, the size and complexity of real-world systems increases the potential number of such alternative optimal solutions exponentially. This is a corollary of the curse of dimensionality: two subnetworks each with alternative solutions may induce a Cartesian product of alternative solutions for the overall network.

There are obvious solutions for the Copper Plate and Network Flow forms of the OPF problem that assign all generation to either Generator 1 or to Generator 2. As a matter of fact, these solutions correspond to optimal solutions for the DC-OPF problem projected to generation variables. When modeled in Pyomo and solved with Gurobi, the DC-OPF problem for this network returns $z^* = 5000$ and $P^* = [100, 0, 0]$. 

There is something that makes even this problem have additional layers of depth: analysis of the tradeoffs in the alternative solutions. The solution to place all of the generation at G1 is optimal on Copper Plate, Network Flow, and DC-OPF, but it will not be optimal for AC-OPF given line losses. This solution will not even be feasible for AC-OPF. Notably, if we had to refine the DC-OPF solution to an AC solution, we would need to take Generator 2 from entirely off to producing at some level. 

There are any number of other secondary concerns that could come up including AC feasibility, market fairness, weather impacts, pollution generation, or availability of fuel. Generating the alternative solutions on OPF problem allows the secondary concerns to be treated outside the optimization model using the finite list of solutions. In its simplest form, this can happen by ranking the alternative solutions on the basis of a given secondary concern; this ranking can be done by additional models/simulators or by a regulatory committee. If both the original model and secondary concerns are linear, then ranking the minimal representation is sufficient to guarantee optimality on the secondary concern.

The minimal representation for the OPF problem will be the extreme points. On the 3-Bus example, all the extreme points for the DC, Network Flow, and Copper Plate OPF problems where generated in Pyomo using the \verb|contrib.enumerate_linear_solutions| method with an optimality gap of 0 \cite{cimor_aos_report}. Since the optimal objective is $z^* = 5000$, $\tau = 5000$ was used. For DC-OPF there are 5 extreme points, 4 for Network Flow, and 2 for Copper Plate. The Copper Plate alternative solutions are the simplest as $P^{(1)} = [100, 0, 0]$ and $P^{(2)} = [0, 100, 0]$. The Network Flow alternative solutions are the same when projected onto the generation variables. The DC-OPF extreme points differ when projected onto just the generation variables, $P^{(1)}$ and $P^{(2)}$ still occur but now $P^{(3)} = [50, 50, 0]$ appears. Since $P^{(3)} \in conv(P^{(1)}, P^{(2)})$, this point implies there is non-equivalence between the DC-OPF alternative solutions versus Network Flow and Copper Plate over the projected out variables.

To see how the extreme points differ on secondary concerns consider the DC to AC feasibility point. The DC-OPF alternative solutions has an exact optimal solution that involves splitting the load between both generators equally. This addresses the point that with any line loss, a single generator active solution will be inherently AC-OPF infeasible with solution adjustment requiring activating an otherwise inactive/cold generator.
\section{Conclusions}
\label{conclusion}
The integration of alternative optimal solutions methods into power systems optimization practice is ongoing. The underlying mathematical optimization modeling tools and solvers continue to advance, and that gives rise to feasible AOS generation methods. For specific applications, MGA-style methods have enhanced ESOM and CapEx modeling results \cite{decarolis_alternatives_mga_energy_futures}. The effort to improve optimal power system results with alternative solutions depends on the operations research and  power system researchers having a common framework on which to build. This paper has endeavored to address the three most important aspects of building that common framework. First and foremost, power systems researchers have to be aware that alternative optimal solutions exist as a concept and given the mathematical tools to apply these methods. The development of the combined sublevel set, minimal representation, and projection framework for analyzing alternative solutions is one clear step in that direction. Similarly, researchers need to be aware of the various methods and software already developed that can be adapted to their specific problems. The literature review herein pulls together such alternative solution references from across power systems, infrastructure planning, agricultural economics, and operations research. Finally, researchers need to understand how alternative solutions look and behave on their problems. It was this last point that motivated the choice of example of Optimal Power Flow as the fundamental example power systems example given use across economic dispatch, unit commitment, transmission switching, capacity expansion, and transmission expansion domains \cite{skolfield_opf_review}. It should now be clear that an optimal solution is not enough on its own and alternative solutions methods can add critical context to optimal power systems applications.






\bibliographystyle{IEEEtran}
\bibliography{IEEEabrv,local_bib}

\begin{thebibliography}{10}
\providecommand{\url}[1]{#1}
\csname url@samestyle\endcsname
\providecommand{\newblock}{\relax}
\providecommand{\bibinfo}[2]{#2}
\providecommand{\BIBentrySTDinterwordspacing}{\spaceskip=0pt\relax}
\providecommand{\BIBentryALTinterwordstretchfactor}{4}
\providecommand{\BIBentryALTinterwordspacing}{\spaceskip=\fontdimen2\font plus
\BIBentryALTinterwordstretchfactor\fontdimen3\font minus \fontdimen4\font\relax}
\providecommand{\BIBforeignlanguage}[2]{{%
\expandafter\ifx\csname l@#1\endcsname\relax
\typeout{** WARNING: IEEEtran.bst: No hyphenation pattern has been}%
\typeout{** loaded for the language `#1'. Using the pattern for}%
\typeout{** the default language instead.}%
\else
\language=\csname l@#1\endcsname
\fi
#2}}
\providecommand{\BIBdecl}{\relax}
\BIBdecl

\bibitem{brill_1979_public_planning_alternatives}
\BIBentryALTinterwordspacing
E.~D. Brill, ``The use of optimization models in public-sector planning,'' \emph{Management Science}, vol.~25, no.~5, pp. 413--422, 1979. [Online]. Available: \url{http://www.jstor.org/stable/2630272}
\BIBentrySTDinterwordspacing

\bibitem{paris_mult_opt_lp}
\BIBentryALTinterwordspacing
Q.~Paris, ``Multiple optimal solutions in linear programming models,'' \emph{American Journal of Agricultural Economics}, vol.~63, no.~4, pp. 724--727, 1981. [Online]. Available: \url{http://www.jstor.org/stable/1241218}
\BIBentrySTDinterwordspacing

\bibitem{VOLL_optimum_is_not_enough}
\BIBentryALTinterwordspacing
P.~Voll, M.~Jennings, M.~Hennen, N.~Shah, and A.~Bardow, ``The optimum is not enough: A near-optimal solution paradigm for energy systems synthesis,'' \emph{Energy}, vol.~82, pp. 446--456, 2015. [Online]. Available: \url{https://doi.org/10.1016/j.energy.2015.01.055}
\BIBentrySTDinterwordspacing

\bibitem{decarolis_alternatives_mga_energy_futures}
\BIBentryALTinterwordspacing
J.~F. DeCarolis, ``Using modeling to generate alternatives (mga) to expand our thinking on energy futures,'' \emph{Energy Economics}, vol.~33, no.~2, pp. 145--152, 2011. [Online]. Available: \url{https://doi.org/10.1016/j.eneco.2010.05.002}
\BIBentrySTDinterwordspacing

\bibitem{Lau_2024_mga_alternatives_energy_applications_review}
\BIBentryALTinterwordspacing
M.~Lau, N.~Patankar, and J.~D. Jenkins, ``Measuring exploration: evaluation of modelling to generate alternatives methods in capacity expansion models,'' \emph{Environmental Research: Energy}, vol.~1, no.~4, p. 045004, 10 2024. [Online]. Available: \url{https://dx.doi.org/10.1088/2753-3751/ad7d10}
\BIBentrySTDinterwordspacing

\bibitem{turan2025oracle_preprint}
\BIBentryALTinterwordspacing
E.~M. Turan, S.~Moret, and A.~Bardow, ``Oracle: A rigorous metric and method to explore all near-optimal designs for energy systems,'' 2025, unpublished. [Online]. Available: \url{https://arxiv.org/abs/2509.26452}
\BIBentrySTDinterwordspacing

\bibitem{paris_economic_foundations_of_symmetric_programming}
Q.~Paris, \emph{Economic foundations of symmetric programming}.\hskip 1em plus 0.5em minus 0.4em\relax Cambridge University Press, 2010.

\bibitem{hpwilliams_model_building}
H.~P. Williams, \emph{Model building in mathematical programming}.\hskip 1em plus 0.5em minus 0.4em\relax Wiley, 2013.

\bibitem{robinson_stability_system_of_inequalities}
\BIBentryALTinterwordspacing
S.~M. Robinson, ``Stability theory for systems of inequalities. part i: Linear systems,'' \emph{SIAM Journal on Numerical Analysis}, vol.~12, no.~5, pp. 754--769, 1975. [Online]. Available: \url{http://www.jstor.org/stable/2156189}
\BIBentrySTDinterwordspacing

\bibitem{robinson_stability_linear_programs}
\BIBentryALTinterwordspacing
------, ``A characterization of stability in linear programming,'' \emph{Operations Research}, vol.~25, no.~3, pp. 435--447, 1977. [Online]. Available: \url{http://www.jstor.org/stable/169931}
\BIBentrySTDinterwordspacing

\bibitem{brill_1982_alternatives_HSJ_hop_skip_jump}
\BIBentryALTinterwordspacing
E.~D. Brill, S.-Y. Chang, and L.~D. Hopkins, ``Modeling to generate alternatives: The hsj approach and an illustration using a problem in land use planning,'' \emph{Management Science}, vol.~28, no.~3, pp. 221--235, 1982. [Online]. Available: \url{http://www.jstor.org/stable/2630877}
\BIBentrySTDinterwordspacing

\bibitem{Convex_Analysis_Rockafellar}
R.~T. Rockafellar, \emph{Convex analysis (Princeton paperbacks)}.\hskip 1em plus 0.5em minus 0.4em\relax Princeton University Press, 1996.

\bibitem{robinson1984_stability_sublevel_set}
\BIBentryALTinterwordspacing
S.~M. Robinson, ``Persistence and continuity of local minimizers,'' \emph{Collaborative Papers}, 1984. [Online]. Available: \url{https://pure.iiasa.ac.at/id/eprint/2571/7/CP-84-005.pdf}
\BIBentrySTDinterwordspacing

\bibitem{royset2021optimization_primer}
J.~O. Royset and R.~J.-B. Wets, \emph{An optimization primer}.\hskip 1em plus 0.5em minus 0.4em\relax Springer, 2021, vol. 440.

\bibitem{convex_opt_theory}
D.~P. Bertsekas, \emph{Convex optimization theory}.\hskip 1em plus 0.5em minus 0.4em\relax Athenea Scientific, 2009.

\bibitem{conforti2014integer}
M.~Conforti, G.~Cornu{\'e}jols, and G.~Zambelli, \emph{Integer Programming}, ser. Graduate Texts in Mathematics.\hskip 1em plus 0.5em minus 0.4em\relax Springer International Publishing, 2014.

\bibitem{lee2000recursive}
\BIBentryALTinterwordspacing
S.~Lee, C.~Phalakornkule, M.~M. Domach, and I.~E. Grossmann, ``Recursive milp model for finding all the alternate optima in lp models for metabolic networks,'' \emph{Computers \& Chemical Engineering}, vol.~24, no. 2-7, pp. 711--716, 2000. [Online]. Available: \url{https://doi.org/10.1016/S0098-1354(00)00323-9}
\BIBentrySTDinterwordspacing

\bibitem{aos_benders_preprint_viens_hart_ferris_2025}
\BIBentryALTinterwordspacing
M.~Viens, W.~E. Hart, and M.~Ferris, ``Extracting alternative solutions from benders decomposition,'' 2025, unpublished. [Online]. Available: \url{https://arxiv.org/abs/2509.08671}
\BIBentrySTDinterwordspacing

\bibitem{Petit2019_alternative_solutions_as_solution_engineering_review_ijoc}
\BIBentryALTinterwordspacing
T.~Petit and A.~C. Trapp, ``Enriching solutions to combinatorial problems via solution engineering,'' \emph{INFORMS Journal on Computing}, vol.~31, no.~3, pp. 429--444, Jul 2019. [Online]. Available: \url{https://doi.org/10.1287/ijoc.2018.0855}
\BIBentrySTDinterwordspacing

\bibitem{Ahanor_Trapp_DiversiTree_IJOC}
\BIBentryALTinterwordspacing
I.~Ahanor, H.~Medal, and A.~C. Trapp, ``{DiversiTree}: A new method to efficiently compute diverse sets of near-optimal solutions to mixed-integer optimization problems,'' \emph{INFORMS Journal on Computing}, vol.~36, no.~1, pp. 61--77, 2024. [Online]. Available: \url{https://doi.org/10.1287/ijoc.2022.0164}
\BIBentrySTDinterwordspacing

\bibitem{danna2007_multiple_solutions_for_mixed_integer_programs_mips_ILOG_CPLEX}
E.~Danna, M.~Fenelon, Z.~Gu, and R.~Wunderling, ``Generating multiple solutions for mixed integer programming problems,'' in \emph{Integer Programming and Combinatorial Optimization}, M.~Fischetti and D.~P. Williamson, Eds.\hskip 1em plus 0.5em minus 0.4em\relax Berlin, Heidelberg: Springer Berlin Heidelberg, 2007, pp. 280--294.

\bibitem{balas_no_good_cuts}
\BIBentryALTinterwordspacing
E.~Balas and R.~Jeroslow, ``Canonical cuts on the unit hypercube,'' \emph{SIAM Journal on Applied Mathematics}, vol.~23, no.~1, pp. 61--69, 1972. [Online]. Available: \url{https://www.jstor.org/stable/2099623}
\BIBentrySTDinterwordspacing

\bibitem{denegre2011interdiction_integer_no_good_cuts_thesis}
S.~DeNegre, ``Interdiction and discrete bilevel linear programming,'' Ph.D. dissertation, Lehigh University, 2011.

\bibitem{PyPSA}
\BIBentryALTinterwordspacing
T.~Brown, J.~H\"orsch, and D.~Schlachtberger, ``{PyPSA: Python for Power System Analysis},'' \emph{Journal of Open Research Software}, vol.~6, no.~4, 2018. [Online]. Available: \url{https://doi.org/10.5334/jors.188}
\BIBentrySTDinterwordspacing

\bibitem{gurobi}
{Gurobi Optimization, LLC}, ``{Gurobi Optimizer Reference Manual},'' 2024.

\bibitem{cplex}
{IBM ILOG CPLEX}, ``{User's Manual for CPLEX},'' 2024.

\bibitem{gams}
\BIBentryALTinterwordspacing
G.~D. Corporation, ``Gams: The general algebraic modeling system,'' 2025, version: 50.1.0. [Online]. Available: \url{https://www.gams.com/}
\BIBentrySTDinterwordspacing

\bibitem{bynum2021pyomo_book}
M.~L. Bynum, G.~A. Hackebeil, W.~E. Hart, C.~D. Laird, B.~L. Nicholson, J.~D. Siirola, J.-P. Watson, and D.~L. Woodruff, \emph{Pyomo--optimization modeling in python}, 3rd~ed.\hskip 1em plus 0.5em minus 0.4em\relax Springer Science \& Business Media, 2021, vol.~67.

\bibitem{cimor_aos_report}
\BIBentryALTinterwordspacing
W.~E. Hart, E.~S. Johnson, C.~A. Phillips, A.~Aquino, B.~Ammari, B.~Arguello, S.~A. Davis, J.~L. Gearhart, C.~D. Laird, C.~L. Mattes \emph{et~al.}, ``{CI-MOR Final Report: Analysis and Validation of Critical Infrastructure Models using Model Order Reduction},'' Sandia National Lab. (SNL-NM), Albuquerque, NM (United States), Tech. Rep., 10 2024. [Online]. Available: \url{https://www.osti.gov/biblio/2480173}
\BIBentrySTDinterwordspacing

\bibitem{frank2012optimal_power_flow_standard_reference_opf}
\BIBentryALTinterwordspacing
S.~Frank, I.~Steponavice, and S.~Rebennack, ``Optimal power flow: A bibliographic survey i: Formulations and deterministic methods,'' \emph{Energy systems}, vol.~3, no.~3, pp. 221--258, 2012. [Online]. Available: \url{https://doi.org/10.1007/s12667-012-0056-y}
\BIBentrySTDinterwordspacing

\bibitem{skolfield_opf_review}
\BIBentryALTinterwordspacing
J.~K. Skolfield and A.~R. Escobedo, ``Operations research in optimal power flow: A guide to recent and emerging methodologies and applications,'' \emph{European Journal of Operational Research}, vol. 300, no.~2, pp. 387--404, 2022. [Online]. Available: \url{https://doi.org/10.1016/j.ejor.2021.10.003}
\BIBentrySTDinterwordspacing

\end{thebibliography}
\appendix
\subsection{OPF Theorem Appendix}
\label{theorem_appendix}
\begin{theorem}
    For $\tau \in \mathR$, the following hold:
    $$\proj_{|P|,|f|}(\levelset(g, X_{DC}, \tau)) \subseteq \levelset(g, X_{NF}, \tau)$$
    $$\proj_{|P|}(\levelset(g, X_{DC}, \tau)) \subseteq \proj_{|P|}(\levelset(g, X_{NF}, \tau)) \subseteq \levelset(g, X_{CP}, \tau)$$
\end{theorem}
\begin{IEEEproof}
    The definitions of $X_{DC}, X_{NF}, X_{CP}$ trivially give:
    $$\proj_{|P|,|f|}(X_{DC}) \subseteq X_{NF}$$
    $$\proj_{|P|}(X_{DC}) \subseteq \proj_{|P|}(X_{NF}) \subseteq X_{CP}$$
    Then by construction of the sublevel set as $g(P) \leq \tau$ applies an identical cut to all the sublevel sets.
\end{IEEEproof}
Since this holds for every fixed value of $\tau$ we immediately get the following corollary:
\begin{corollary}
    For $z^* = \min_{(P,f,\theta)} g(P)$, then both:
    \begin{align*}
        \proj_{|P|,|f|}(\levelset(g, X_{DC}, z^*)) \subseteq & \  \levelset(g, X_{NF}, z^*),
    \end{align*}
    and
    \begin{align*}
        \proj_{|P|}(\levelset(g, X_{DC}, z^*)) &\subseteq \  \proj_{|P|}(\levelset(g, X_{NF}, z^*)) \\
        & \subseteq \levelset(g, X_{CP}, z^*)
    \end{align*}
    hold.
\end{corollary}

\end{document}